\documentstyle[12pt]{article}

\newtheorem{lemma}{Lemma}
[section]

\newtheorem{thm}{Theorem}
[section]
\newtheorem{cor}{Corollary}
[section]

\newcommand{\mysection}[1]{\section{#1}\setcounter{equation}{0}}

\newcommand{\la}{\lambda}

\newcommand{\ptl}{\partial}

\newcommand{\ra}{\rightarrow}
\newcommand{\lra}{\longrightarrow}

\def\a{\alpha}

\def\R2{{\Bbb R}^2}
\def\Bbb{\mathbb}

\def\beq{\begin{equation}}
\def\eeq{\end{equation}}
\def\ba{\begin{array}}
\def\ea{\end{array}}
\def\barr{\begin{eqnarray}}
\def\earr{\end{eqnarray}}

\def\part{\partial}
\def\fr{\frac}

\newcommand{\ls}{\setlength{\baselineskip}{18pt}
                      \setlength{\parskip}{3mm} }

\title {The H\"older continuity of a class of 3-dimension ultraparabolic equations}
\author{ WANG Wendong and ZHANG Liqun
\thanks{ The research is partially supported by the
Chinese NSF under grant 10325104. Email: wangwendong@amss.ac.cn and
lqzhang@math.ac.cn} }

\date{Institute of Mathematics, AMSS, Academia Sinica,
Beijing\\  }

\begin{document}

\maketitle

\begin{abstract}
We obtained the $C^{\a}$ continuity for weak solutions of a class of
ultraparabolic equations with measurable coefficients
 of the form $${\ptl_t \,u}= \ptl_x(a(x,y,t)\ptl_x \,u
)+b_0(x,y,t)\ptl_x u+b(x,y,t)\ptl_y u,$$ which generalized our
recent results on KFP equations.
\end{abstract}

{\small keywords: Hypoelliptic, ultraparabolic equations, H\"older
regularity }

\pagenumbering{arabic}
\mysection{Introduction} \ls \noindent

Consider a class of ultraparabolic operator on ${R}^{2+1}$:
$$
\displaystyle  Lu \equiv \ptl_x(a(x,y,t)\ptl_x \,u
)+b_0(x,y,t)\ptl_x u+b(x,y,t)\ptl_y u-{\ptl_t \,u}=0,\leqno(1.1)
$$
where $(x,y,t)=z\in \Omega\subset{ R}^{2+1}$,  $a(z)$, $b_0(z)$ and
$b(z)$ is real, measurable functions. We assume that $b(z)$ is twice
differentiable, and there exists a positive constant $\mu$ such that
for $z\in \Omega$,
$$\mu<a(z)<\mu^{-1},\quad \quad\fr{\ptl b(z)}{\ptl x} \neq 0,
\quad  \quad |b|_{C^2}+|b_0|_{\infty}\leq \mu^{-1}. \leqno (1.2)$$
Also, we denote $$ L_0u = \ptl_x^2 \,u +x\ptl_y u-{\ptl_t
\,u}=0,\leqno(1.3)$$
$$ L_1u = \ptl_x(a(x,y,t)\ptl_x \,u)
 +x\ptl_y u-{\ptl_t \,u}=0,\leqno(1.4)$$
 and
$$
L_2 u = \ptl_x(a(x,y,t)\ptl_x \,u )+b_0(x,y,t)\ptl_x u+x\ptl_y
u-{\ptl_t \,u}=0.\leqno(1.5)
$$

We remark that the equation (1.3) and (1.4) are the examples of
3-dimension homogeneous Kolmogorov-Fokker-Planck equations (or KFP
equations). The condition $\fr{\ptl b}{\ptl x} \neq 0$ ensures (1.1)
satisfies H\"ormander's hypoellipticity conditions,
$$
{\rm rank \,\, Lie}(\ptl_x, b\ptl_y-\ptl_t)(z)=3, \quad \forall z\in
\Omega .
$$
The study of regularity of the KFP equation has a long history, and
the earlier works are mainly on the Schauder type estimates. The
study of regularity of weak solutions is begun in recent years. A
recent paper of Pascucci and Polidoro [6], has proved that the Moser
iterative method still works for the class of KFP equations with
measurable coefficients. By the same technique, Cinti, Pascucci,
Polidoro [1] consider a class of nonhomogeneous KFP equations,  and
Cinti, Polidoro [2] deal with a more general ultraparabolic
equation. Their results show that for a non-negative sub-solution
$u$ of the ultraparabolic equation, $L^{\infty}$ norm of $u$ is
bounded by the $L^p$ norm ($p \ge 1$). The second author [10], [11]
has proved $C^{\a}$ property of weak solutions by Kruzhkov's
approach for homogeneous KFP equations, and the authors deal with
nonhomogeneous KFP equations in [7]. By simplifying the cut-off
function and generalizing their earlier arguments, the authors [8]
have considered more general ultraparabolic equations whose
fundamental solution is implicit. We are not try to review the
detailed history, but focus on the study of the H\"older continuity
of a simple looking case. In this paper, we give another
generalization of KFP equations in $R^{2+1}$ and consider the
hypoelliptic operator as L in (1.1).

We say that $u$ is a $weak \,solution$ if it satisfies (1.1) in the
distribution sense, that is for any $\phi \in C^{\infty}_0(\Omega)$,
$\Omega$ is a open subset of $R^{2+1}$, then
$$
\int_{\Omega} \phi (b_0\ptl_x+b\ptl_y-\ptl_t) u- a\ptl_x\phi\ptl_x u
= 0, \leqno(1.6)
$$
and $u$, $\ptl_xu$, $b \ptl_y-\ptl_t u \in L^2_{\rm loc}(\Omega).$

Our main result is the following theorem:
\begin{thm}
Under the assumption (1.2), the weak solution of (1.1) is H\"older
continuous.
\end{thm}

\mysection{Some Preliminary and Known Results}

We follow the earlier notations to give some basic known properties
related to our problems. For more details of the subject, we refer
to Pascucci and Polidoro [6] and Lanconelli and Polidoro [5].

Let $ B= \left(
\begin{array}{cc}
0& 1  \\
0 & 0
\end{array}
\right) $, and $E(\tau)=\rm{exp}(-\tau B^T)=$ $ \left(
\begin{array}{cc}
1& 0  \\
-\tau & 1
\end{array}
\right) $.

For $(x,y,t), (\xi,\eta,\tau) \in R^{2+1}$, set
$$(x,y,t)\circ (\xi,\eta,\tau)=((\xi,\eta)+E(\tau)(^x_y),t+\tau),$$
then $(R^{2+1}, \circ)$
is a Lie group with identity element $(0,0)$, and the inverse of an
element is $(x,y,t)^{-1}=(-E(-t)(^x_y),-t)$. The left
translation by $(\xi,\eta,\tau)$ given by
$$(x,y,t)\mapsto (\xi,\eta,\tau)\circ (x,y,t),$$
is a invariant translation to the operator $L_0$. The associated
dilation to operator $L_0$  is given by
$$
\delta_t=diag (t ,t^3 ,t^2),
$$
where $t$ is a positive parameter,
and the homogeneous dimension of
$(R^{2+1},\circ)$ with respect to the dilation $\delta_t$ is 6.

The norm in $R^{2+1}$, related to the group of translations and
dilation to the equation is defined by $||(x,y,t)||=r,$ if $r$ is
the unique positive solution to the equation
$$
\fr{x^2}{r^{2}}+\fr{y^2}{r^{6}} +\fr{t^2}{r^4}=1,
$$
where $(x,y,t) \in R^{2+1}\setminus \{0\}$ . And $||(0,0)||=0$.
Obviously $$\|\delta_{\mu}(x,y,t)\|=\mu\|(x,y,t)\|,$$ for all $(x,y,t)\in
R^{2+1}$.

The ball centered at a point $(x_0,t_0)$ is defined by
$${\cal B}_r(x_0,t_0)=\{(x,t)|\quad ||(x_0,t_0)^{-1}\circ (x,t)||\leq r\},$$
and
$${\cal B}^-_r(x_0,t_0)={\cal B}_r(x_0,t_0)\cap\{t<t_0\}.$$
For convenience, we sometimes use the cube instead of the balls. The
cube at point $(0,0)$ is given by
$$
{\cal C}_r(0,0)=\{(x,y,t)|\quad |x|\leq r, |y|\leq 8r^{3},\quad
|t|\leq r^2\}.
$$
It is easy to see that there exists a constant $\Lambda$ such that
$$
{\cal C}_{\fr r {\Lambda}}(0,0)\subset{\cal B}_r(0,0)\subset{\cal
C}_{\Lambda r}(0,0).
$$

For $ L_0$, the fundamental solution $\Gamma(\cdot,\zeta)$ with pole
in $\zeta=(\xi,\eta,\tau)\in R^{2+1}$ is smooth except at the
diagonal of $R^{2+1}\times R^{2+1}$. It has the following form at
$\zeta=0$,
$$
\Gamma(z)=\Gamma(z,0)= \Big {\{}
\begin{array}{cc}{\frac{\sqrt 3}{2\pi t^2}\exp[-\fr1 t
(x^2+\fr3txy+\fr{3}{t^2}y^2)]} & {\rm if} \quad t>0,
\\
0 & {\rm if} \quad t\leq 0; \end{array}\leqno(2.1)
$$
and
$$
\Gamma(z,\zeta)= \Big {\{}
\begin{array}{cc}{\frac{\sqrt 3}{2\pi (t-\tau)^2}
\exp[-\fr{x^2+x\xi+\xi^2}{ t-\tau}-\fr{3(x+\xi)(y-\eta)}{(t-\tau)^2}
-\fr{3(y-\eta)^2}{(t-\tau)^3}]} & {\rm if} \quad t>\tau,
\\
0 & {\rm if} \quad t\leq \tau. \end{array}\leqno(2.2)
$$
Obviously we can derive from the above formula,
$$
\hspace{2pt}\int_{R^2}\Gamma(x,y,t;\xi,\eta,\tau)dx dy=
\int_{R^2}\Gamma(x,y,t;\xi,\eta,\tau)d\xi d\eta=1,\quad{\rm if}\quad
t>\tau,\leqno(2.3)$$ and
$$\hspace{2pt}\Gamma(\delta_{\mu}\circ z)=\mu^{-4}\Gamma(z),\quad
\forall z\neq 0, \,\mu>0.\leqno(2.4)$$

A $weak$ $sub$-$solution$ of (1.5) in a domain $\Omega$ is a
function $u$ such that $u$, $\ptl_x u$, $(x\ptl_y-\ptl_t)u \in
L^2_{loc}(\Omega)$ and for any $\phi \in C^{\infty}_0(\Omega)$,
$\phi \geq 0$,
$$
\int_{\Omega} \phi (b_0\ptl_x+x\ptl_y-\ptl_t)u-(\ptl_x u) a
\ptl_x\phi \geq 0. \leqno(2.5)
$$

We recall a result of Pascucci and Polidoro obtained by using the
Moser's iterative method (see [6]) states as following
\begin{lemma}
Let $u$ be a non-negative weak sub-solution of (1.4) in $\Omega$.
Let $(x_0,t_0)\in \Omega$ and $\overline{{\cal
B}^-_r(x_0,t_0)}\subset \Omega$ and $p \geq 1$. Then there exists a
positive constant $C$ which depends only on the operator $L$  such
that, for $0 < r\leq 1$
$$
\sup_{{\cal B}^-_{\fr r 2}(x_0,t_0)} u^p \leq \fr
{C}{r^{6}}\int_{{\cal B}^-_r(x_0,t_0)} u^p,\leqno(2.6)
$$
provided that the last integral converges.
\end{lemma}
The second author [10] proved the following result.
\begin{thm}
If u is a weak solution of (1.4), then u is $H\ddot{o}lder$
continuous.
\end{thm}

Using the same technique, we can obtain the similar result to the
equation (1.5).
\begin{thm}
If  u is a weak solution of (1.5), then u is $H\ddot{o}lder$
continuous.
\end{thm}
In section 3, we shall sketch the proof of this theorem first. We
mainly focus on the proof of the oscillation estimates. Then we give
a transformation as Weber in [9] and complete the proof of Theorem
1.1.

We make use of a classical potential estimates (see (1.11) in [3])
here to prove the ${\rm Poincar\acute{e}}$ type inequality.

\begin{lemma}
Let $(R^{N+1},\circ)$ is a homogeneous Lie group of homogeneous
dimension $Q+2$, $\a \in (0, Q+2)$ and $G \in C(R^{N+1}\setminus
\{0\})$ be a $\delta_{\mu}$-homogeneous function of degree $\a-Q-2$.
If $f \in L^p(R^{N+1})$ for some $p \in (1,\infty)$, then
$$
G_f(z)\equiv \int_{R^{N+1}} G(\zeta ^{-1}\circ z)f(\zeta)d\zeta,
$$
is defined almost everywhere and there exists a constant
$C=C(Q,p)$ such that
$$
||G_f||_{L^q(R^{N+1})}\leq C \max_{||z||=1} |G(z)|\quad
||f||_{L^p(R^{N+1})},\leqno(2.7)
$$
where $q$ is defined by
$$
\fr 1q =\fr 1p-\fr{\a}{Q+2}.
$$
\end{lemma}
\begin{cor} Let $f\in L^2(R^{2+1})$,  and recall the definitions in [6]
$$
\Gamma(f)(z)=\int_{R^{2+1}}\Gamma(z,\zeta)f(\zeta) d\zeta, \qquad
\forall z\in R^{2+1},
$$
and
$$
\Gamma(\ptl_{\xi}f)(z)=-\int_{R^{2+1}}\ptl_{\xi}\Gamma(z,\zeta)f(\zeta)
d\zeta, \qquad \forall z\in R^{2+1},
$$
 then exists an absolute constant $C$
such that
$$
\|\Gamma(f)\|_{L^{2\tilde{k}}(R^{2+1})}\leq C\|f\|_{L^2(R^{2+1})},
\leqno(2.8)
$$
and
$$
\|\Gamma(\ptl_{\xi}f)\|_{L^{2k}(R^{2+1})}\leq C\|f\|_{L^2(R^{2+1})},
\leqno(2.9)
$$
where $\tilde{k}=3$, $k=\fr{3}{2}$.
\end{cor}

\mysection{Proof of Theorem 2.2}

Outline of the proof of Theorem 2.2:

Step 1: $L^{\infty}$ estimate via Moser iteration. It can be checked
that the same Caccioppoli type inequality holds ( See Theorem 3.1,
[6]), since
$$\int_{B_1}\psi^2b_0\ptl_x(v^2)\leq \fr12
\int_{B_1}\psi^2av^2+C(\mu)\int_{B_1}\psi^2v^2.$$

In order to use the Moser iteration, one need to prove a Sobolev
type inequality. It can be proved that the Sobolev type inequality
holds for non-negative weak sub-solution ( See Theorem 3.3 in [6]).
Here one may deal with $\int_{B_r}[\Gamma(z,\cdot)v b_0\ptl_x
\psi](\zeta)$ as $I_2$ in [6]. Then one can obtain the $L^{\infty}$
estimate as in Lemma 2.1.

Step 2: Oscillation estimates.

This is obtained in Lemma 3.6. We shall focus on this parts in the
following discussions.

Step 3: H\"older regularity.

This is followed by the oscillation estimated by a standard
argument.

Now we turn to the proof of main results. We may consider the local
estimate at a ball centered at $(0,0)$, since the equation (1.4) is
invariant under the left translation when $a$ is constant. We follow
the same route as [10], [7] and [8]. For convenience, we consider
the estimates in the following cube, instead of ${\cal B}^-_r$,
$$\ba{lll}
{\cal C}_r^{-}=\{(x,y,t)| &-r^2\leq t < 0,\, |x|\leq r,\, |y|\leq (2
r)^3\}.\ea
$$
Let
$$
K_r=\{(x,y)| \,\, |x|\leq r ,\, |y|\leq ( 2r)^3\}.
$$
Let $0<\a, \beta<1$ be constants, for fixed $t$ and $h$, we denote
$$
{\cal N}_{t,h}=\{(x,y)|\,(x,y)\in K_{\beta r} ,\, u(x,y,t) \geq h\}.
$$
 We sometimes abuse the notations of ${\cal
B}^-_r$ and ${\cal C}_r^-$, since there are equivalent.

\begin{lemma}
Suppose that $u(x,t)\geq 0$ be a solution of equation (1.5) in
${\cal C}^-_r$ centered at $(0,0)$ and
$$
mes\{(x,t)\in {\cal C}^-_r, \quad u \geq 1\} \geq \fr 1 2 mes ({\cal
C}^-_r),
$$
then there exist constants $\a$, $\beta$ and $h_1$, $0<\a, \beta, h_1<1$,
where $h_1$ only depends on $\mu$, such that for almost all
$t\in (-\a r^2,0)$ and $0<h<h_1$
$$
mes\{{\cal N}_{t,h}\} \geq \fr {1}{11}mes\{ K_{\beta r}\}
.
$$
\end{lemma}
{\it Proof:} Let
$$
v=\ln^+(\fr{1}{u+h^{\fr 9 8}}),
$$
where $h$ is a constant, $0<h<1$, to be determined later. Then $v$
at points where $v$ is positive, satisfies
$$
\displaystyle \ptl_x(a(x,y,t)\ptl_x \,v )-a(\ptl_x
v)^2+b_0(x,y,t)\ptl_x v+x\ptl_y v-{\ptl_t \,v}=0.\leqno(3.1)
$$
Let $\eta(s)$ be a smooth cut-off function so that
$$
\eta(s)=1,\quad \hbox {for} \quad s< \beta r,
$$
$$
\eta(s)=0,\quad \hbox {for} \quad s\geq r.
$$
Moreover, $0\leq\eta \leq 1$ and $|\eta'|\leq \fr {2}{(1-\beta)r}$.

Multiplying $\eta(|x|)^2$ to (3.1) and integrating by parts on
$K_r\times(\tau,t)$
$$
\ba{lllllllll} &&\int_{K_{\beta r}} v(x,y,t)dx dy
+\int_\tau^t \int_{K_{ r}}\eta^2 \,a
|\ptl_xv|^2d x dyds\\ \\
&\leq& \int_\tau^t \int_{K_{ r}}\eta^2 (\ptl_x(a(x,y,t)\ptl_x \,v
)+b_0(x,y,t)\ptl_x v+x\ptl_y v)dxdyds\\\\&&+\int_{K_{r}}
v(x,y,\tau)dx dy\\\\&\leq& \fr {C(\mu)}
{\beta^4(1-\beta)^2}|K_{\beta r}|+\int_\tau^t
\int_{K_{r}}(\fr12\eta^2 \,a |\ptl_xv|^2+x\ptl_y v\eta^2)dx dy ds
\\\\&&+\int_{K_{r}} v(x,y,\tau)dx dy,\qquad a.e. \quad\tau, t\in(-r^2,0).\ea\leqno(3.2)
$$
Then
$$
\ba{llllll}  |\int_{K_{r}}x\ptl_y v \eta^2|
&=&|\int_{|x|\leq r}xv\eta^2\mid_{y=-8\beta^3r^3}^{8\beta^3r^3}dx|\\\\
&\leq& \fr14r^{-2}\beta^{-4} |K_{\beta r}|\ln(h^{-\fr 9
8}). \ea$$ Integrating by t to $I_B$, we have
$$
\ba{llllll}|\int_\tau^t\int_{K_{r}}x\ptl_y v
\eta^2|  \leq \fr14\beta^{-4} |K_{\beta
r}|\ln(h^{-\fr 9 8}). \ea\leqno (3.3)
$$
We shall estimate the measure of the set ${\cal N}_{t,h}$. Let
$$
\nu(t)=mes\{(x,y)|\quad (x,y)\in K_r, \, u(x,y,t)\geq
1\}.
$$
By our assumption, for $0<\a< \fr 12$
$$
\fr 12 r^2 mes(K_{r})\leq \int_{-r^2}^0
\nu(t)dt=\int_{-r^2}^{-\a r^2}\nu(t)dt+\int_{-\a r^2}^{0}\nu(t)dt,
$$
that is
$$
\int_{-r^2}^{-\a r^2}\nu(t)dt\geq (\fr 12-\a)r^2
mes(K_{r}),
$$
then there exists a $\tau \in (-r^2,-\a r^2)$, such that
$$
\mu(\tau)\geq (\fr 12-\a)(1-\a)^{-1}
mes(K_{r}).
$$
By noticing $v=0$ when $u\geq 1,$ we have
$$
\int_{K_{r}} v(x,y,\tau)dx dy\leq \fr
12(1-\a)^{-1}mes(K_{r})\ln(h^{-\fr 9 8}).\leqno(3.4)
$$
Now we choose $\a$ (near zero), and $\beta$ (near one),  such that
$$
\fr1{4\beta^4}+\fr{1}{2\beta ^{4}(1-\a)}\leq \fr 4 5,\leqno(3.5)
$$
and fix them from now on.

By (3.2), (3.3), (3.4) and (3.5), we deduce
$$\ba{lll}
\int_{K_{\beta r}} v(x,y,t)d x dy\\ \\
\leq [\fr {C(\mu)} {\beta^4(1-\beta)^2} +\fr 45\ln(h^{-\fr 9
8})]mes(K_{\beta r}).\ea\leqno(3.6)
$$
When $(x, y)\notin {\cal N}_{t,h},$, we have
$$\ln(\fr 1 {2h})\leq \ln^+(\fr{1}{h+h^{\fr 9 8}})\leq v,$$
then
$$\ln(\fr 1
{2h})mes(K_{\beta r}\setminus {\cal N}_{t,h})\leq
\int_{K_{\beta r}} v(x,y,t)d x dy.$$ Since
$$
\fr{C+{\fr 45}\ln(h^{-\fr 98})}{\ln(h^{-1})}\lra \fr
9{10},\qquad\hbox{as} \quad h\ra 0,
$$
then there exists constant $h_1$ such that for $0<h<h_1$ and $t
\in(-\a r^2,0)$
$$
mes(K_{\beta r}\setminus {\cal N}_{t,h})\leq \fr
{10}{11}mes(K_{\beta r}).
$$
Then we proved our lemma.

Let $\chi(s)$ be a $C^{\infty}$ smooth function given by
$$\ba{ll}
\chi(s)=1 \qquad if \quad s\leq {\theta^{\fr 1 {6}}} r,\\
\chi(s)=0 \qquad if \quad s> r, \ea
$$
where $\theta>0$ is a constant, to be determined in Lemma 3.4, and
${\theta^{\fr 1 {6}}}<\fr {1}{2}$. Moreover, we assume that
$$
0\leq -\chi'(s) \leq \fr{2}{(1 -{\theta^{\fr 1 {6}}})r},\quad
|\chi''(s)|\leq \fr{C}{r^2},
$$
and for any $\beta_1, \beta_2,$ with $\theta^{\fr 1
{6}}<\beta_1<\beta_2<1,$ we have $$|\chi'(s)|\geq
C(\beta_1,\beta_2)r^{-1}>0,$$ if $\beta_1r\leq s\leq \beta_2r.$

For $(x,y)\in R^2,$ $t\leq 0$, we set $${\mathcal Q} =\{(x,y,t)|
-r^2\leq t\leq 0,\, |x|\leq \fr r {\theta}, \,|y|\leq {\fr {r^3}
{\theta}}\}.$$ We define the cut off functions by
$$
\phi_0(x,y,t)=\chi([\theta^2y^2-6 t r^{4} ]^{\fr {1} {6}}),
$$
$$
\phi_1(x,y,t)=\chi(\theta |x|),
$$
$$
\phi(x,y,t)=\phi_0 \phi_1. \leqno(3.7)
$$

\begin{lemma}
By the definition of $\phi$ and the above arguments, we have
 $$(x\ptl_y-\ptl_t)\phi_0(z)\leq 0, \quad \rm{for}\quad z\in {\mathcal Q}.\leqno(3.8)$$
And since $\theta^{\fr16}<\fr1{2}$, we have\\
(1) $\phi(z)\equiv 1,$ in ${\cal B}^-_{\theta r}$,\\
(2) $\rm{supp}\phi\bigcap \{(x,y,t)|(x,y)\in R^2,t\leq 0\}\subset {\mathcal Q}$,\\
(3) there exists a constant $\a_1$, $0<\a_1<\min \{\a,\fr 1{12}\}$, such that
$$\{(x,y,t)| -\a_1r^2\leq t < 0, (x,y)\in K_{\beta r}
\}\subseteq \rm{supp}\phi, $$moreover, $0<\phi_0(z)<1,$ for $z\in
\{(x,y,t)| -\a_1r^2\leq t \leq -\theta r^2, (x,y)\in K_{\beta r}
\}$.
\end{lemma}
{\it Proof:} By the definition of $\phi_0$, we attain
$$\ba{llllllllll} (x\ptl_y-\ptl_t) \phi_0 &=
\chi'([\theta^2 y^2-6 tr^4]^{\fr{1}{6}})\fr{1}{6}[\theta^2 y^2-6 tr^4]^{-\fr{5}{6}} [6 r^{4}
+2\theta^2 xy]\leq 0.
 \ea
$$
When $\theta<\fr16$, we can check that obviously (1) holds. We
notice that either $|x|\geq \fr r {\theta}$, or $|y|\geq \fr
{r^3}{\theta}$, or $t\leq -r^2$, then $\phi$ vanishes, hence we
obtain (2). When $(x,y)\in K_{\beta r}$, then $\phi_1>0$ and we can
choose $\theta<\fr1{64}$ and $t$ small, for example, $t>-\a_1r^2$,
such that $\theta^2 y^2-6 tr^4<r^6$, then we obtain (3).

 Now we have the following ${\rm
Poincar\acute{e}}$'s type inequality.
\begin{lemma}
Let $w$ be a non-negative weak sub-solution of (1.5) in ${\cal
B}_1^-$. Then there exists an absolute constant $C$, such that for
$r<\theta<1$
$$
\int_{{\cal B}^-_{\theta r}}(w(z)-I_0)_+^2\leq C\theta^2
r^2\int_{{\cal B}^-_{\fr r {\theta}}}|\ptl_x w|^2, \leqno(3.9)
$$
where
$
I_0=sup_{{\cal B}^-_{\theta r}}I_1(z),
$
and
$$
I_1(z)=\int_{{\cal B}^-_{\fr r {\theta}}}
[-\Gamma(z,\zeta)w(\zeta)(\xi\ptl_{\eta}-\ptl_{\tau})\phi(\zeta)-
\ptl_{\xi}^2\phi(\zeta)\Gamma(z,\zeta) w(\zeta)]d\zeta,\leqno(3.10)
$$
where $\Gamma$ is the fundamental solution of $L_0$,  and $\phi$ is
given by (3.7).
\end{lemma}
 {\it Proof:} We represent $w$ in terms of the fundamental
solution of $\Gamma$, i.e.
$$\varphi(z)=-\int_{R^{2+1}}\Gamma(z,\zeta)L_0\varphi(\zeta)d\zeta,
 \quad \forall\varphi\in C_0^{\infty}(R^{2+1}).$$ By an approximation of
 $\phi$ and integrating by parts,  for $z \in {\cal B}^-_{\theta r}$,
we have
$$\ba{llll}
w(z)&=\int_{{\cal B}^-_{\fr r {\theta}}}  [\langle
\ptl_{\xi}(w\phi)(\zeta),\ptl_{\xi}\Gamma(z,\zeta)\rangle
-\Gamma(z,\zeta)(\xi\ptl_{\eta}-\ptl_{\tau})(w\phi)(\zeta)]d\zeta \\
\\&= I_1(z)+I_2(z)+I_3(z),\ea\leqno(3.11)
$$
where
$$
I_1(z)=\int_{{\cal B}^-_{\fr r {\theta}}}
[-\Gamma(z,\zeta)w(\xi\ptl_{\eta}-\ptl_{\tau})\phi+\langle\ptl_{\xi}{\phi},\ptl_{\xi}\Gamma(z,\zeta)\rangle
w+\Gamma(z,\zeta)\langle\ptl_{\xi}w,\ptl_{\xi}\phi\rangle]d\zeta,$$

$$ \ba{lllllll}
I_2(z)&=&\int_{{\cal B}^-_{\fr r {\theta}}} [\langle
(1-a)\ptl_{\xi}w,\ptl_{\xi}\Gamma(z,\zeta)\rangle\phi-\Gamma(z,\zeta)\langle
(a+1)\ptl_{\xi}w,\ptl_{\xi}\phi\rangle\\\\& &+\Gamma(z,\zeta)\phi
b_0\ptl_{\xi}w]d\zeta =I_{21}+I_{22}+I_{23},  \ea
$$
and $$ I_3(z)=\int_{{\cal B}^-_{\fr r {\theta}}}  [\langle a
\ptl_{\xi}w,\ptl_{\xi}(\Gamma(z,\zeta)\phi)\rangle
-\Gamma(z,\zeta)\phi(b_0\ptl_{\xi}+\xi\ptl_{\eta}-\ptl_{\tau})w]d\zeta
.
$$
Note that ${\rm supp}\phi\bigcap\{\tau\leq 0\}\subset {\mathcal
Q}\subset \overline{{\cal B}^-_{\fr r {\theta}}}$, $z \in {\cal
B}^-_{\theta r}$ and $\langle
\ptl_{\xi}{\phi},\ptl_{\xi}\Gamma(z,\zeta)\rangle $ vanishes in a
small neighborhood of $z$. Integrating by parts we obtain $I_1(z)$
as in (3.10).

From our
assumption, $w$ is a weak sub-solution of (1.5), and $\phi$ is a
test function of this semi-cylinder. In fact, we let
$$
\tilde{\chi}(\tau)=\left\{
\begin{array}{lll} 1\quad &\tau\leq
0,\\1-n\tau\quad &0\leq \tau \leq 1/n,\\ 0\quad &\tau\geq
1/n.\end{array}\right.$$ Then
$\tilde{\chi}(\tau)\phi\Gamma(z,\zeta)$ can be a test function (see
[6]). As $n\rightarrow \infty$, we obtain $\phi\Gamma(z,\zeta)$ as a
legitimate test function, and $I_3(z)\leq 0$. Then in ${\cal
B}^-_{\theta r}$,
$$
0\leq (w(z)-I_0)_+\leq I_2(z).
$$
By Corollary 2.1 we have
$$
||I_{21}||_{L^2({\cal B}^-_{\theta r})}\leq C\theta
r||I_{21}||_{L^{3}({\cal B}^-_{\theta r})}\leq C \theta
r||\ptl_{\xi}w||_{L^2({\cal B}^-_{{\fr r {\theta}}})}.\leqno(3.12)
$$
Similarly
$$||I_{2i}||_{L^2({\cal B}^-_{\theta r})}\leq C \theta^2
r||\ptl_{\xi}w||_{L^2({\cal B}^-_{{\fr r {\theta}}})},\quad
i=2,3$$then we proved our lemma.

Now we apply Lemma 3.3 to the function $ w= \ln^+\fr{h}{u+h^{\fr
98}}. $ If $u$ is a weak solution of (1.5), then $w$ is a weak
sub-solution. We estimate the value of $I_0$.
\begin{lemma}
Under the assumptions of Lemma 3.3, there exist constants $\theta$,
$\lambda_0$,  $\la_0<1$ only depends on constants $\a$ and $\beta$,
 such that for $r<\theta$
$$
|I_0|\leq \lambda_0 \ln(h^{-\fr 1 8}).\leqno(3.13)
$$
\end{lemma}
{\it Proof:} We first come to estimate the second term of $I_1(z)$ and as before,
denote $z=(x,y,t)$ and $\zeta=(\xi,\eta,\tau)$. Note
 $z \in
{\cal B}^-_{\theta r}$, we have
$$
\ba{llllllllllll}
&& \int_{{\cal B}^-_{\fr r {\theta}}}
[|\ptl_{\xi}^2\phi(\zeta)|\Gamma(z,\zeta) w]d\zeta\\\\
&\leq& r^2sup_{\xi\in {\rm
supp}(\ptl_{\xi}\phi)}|\ptl_{\xi}^2\phi|(\zeta)\ln (h^{-\fr 1
8}).\quad (By\,\,\,(2.3))\ea$$ We only need to estimate
$|\ptl_{\xi}^2\phi_1|$. Since
$$|\ptl_{\xi}\phi_1|=|\theta\chi'(\theta
|\xi|)\partial_{\xi}|\xi||\leq 4\theta r^{-1},\quad
|\ptl_{\xi}^2\phi_1|\leq C\theta^{2}r^{-2}. $$ Hence
$$| \int_{{\cal B}^-_{\fr r {\theta}}}[-\ptl_{\xi}^2\phi\Gamma(z,\zeta) w]d\zeta |
\leq C_3\theta^2 \ln (h^{-\fr 1 8})
$$
where $C_3$ is an absolute constant.

Now we let $w\equiv 1$, then for $z \in {\cal B}^-_{\theta r}$,
(3.11) gives
$$
\ba{llll} 1=\int_{{\cal B}^-_{\fr r {\theta}}}[ -
\phi_1\Gamma(z,\zeta)(\xi\ptl_{\eta}-\ptl_{\tau})\phi_0]d\zeta +
\int_{{\cal B}^-_{\fr r {\theta}}}[-\ptl_{\xi}^2\phi(\zeta)\Gamma(z,\zeta)]d\zeta.\ea\
\leqno(3.14)
$$
By (3.8) in Lemma 3.2, we know that
$$
-\phi_1\Gamma(z,\zeta)(\xi\ptl_{\eta}-\ptl_{\tau})\phi_0 \geq 0.
$$
We only need to prove
$-\phi_1\Gamma(z,\zeta)(\xi\ptl_{\eta}-\ptl_{\tau})\phi_0$ has a
positive lower bound in a domain which $w$ vanishes, and this bound
independent of $\,r$ and small $\theta$. So we can find a $\la_0,$
$0<\la_0<1$, such that this lemma holds.

For $z\in B^-_{\theta r}$, set
$$\zeta\in Z=\{(\xi,\eta,\tau)| -\a_1 r^2\leq \tau\leq -\fr {\a_1}
{2}r^2,\,(\xi,\eta)\in K_{\beta r},
\,w(\xi,\eta,\tau)=0\},\leqno(3.15)$$ then  by Lemma 3.1,
$|Z|=C(\a_1, \beta)r^{6}$.

We note that when $\zeta=(\xi,\eta,\tau)\in Z $ and
$\theta<\fr1{64}$, $w(\zeta)=0,$ $\phi_1(\zeta)=1$, and
$$|\chi'([\theta^2y^2-6 t r^{4}]^{\fr{1}{6}})|\geq C(\a_1)r^{-1}>0.$$
Consequently
$$
\ba{llllllllllll} \int_Z [-\phi_1\Gamma(z,\zeta)(\xi\ptl_{\eta}-\ptl_{\tau})\phi_0]\,d \zeta \\
\\= -\int_Z \phi_1\Gamma(z,\zeta)
\chi'([\theta^2\eta^2-6 \tau r^{4}]^{\fr{1}{6}})\fr{1}{6}[\theta^2\eta^2-6 \tau r^{4}]^{-\fr{5}{6}} [6 r^{4}
+2\theta^2\xi\eta]d\zeta\\
\\ \geq C(\a_1)   \int_Z
r^{-2}\Gamma(\zeta^{-1}\circ z,0)d\zeta
=C(\a,\beta)=C_4>0, \ea
$$
where we have used $\Gamma(z,\zeta)\geq Cr^{-4}$, as $\tau \leq
-\fr{\a_1}{2}r^2$ and $z\in B^-_{\theta r}$. In fact, by (2.2) one
can obtain this result easily.

By (3.14) and  $$I_0(z)=\sup_{{\cal B}^-_{ {\theta r}}}\int_{{\cal
B}^-_{\fr r {\theta}}\setminus Z}
[-\Gamma(z,\zeta)w(\zeta)(\xi\ptl_{\eta}-\ptl_{\tau})\phi(\zeta)-\ptl_{\xi}^2\phi(\zeta)\Gamma(z,\zeta)
w(\zeta)]d\zeta, $$ we have
$$
|I_0|\leq (1-C_4+C_3\theta^{2})\ln(h^{-{\fr 1 8
}})+C_3\theta^{2}\ln(h^{-{\fr 1 8 }}) .
\leqno(3.16)
$$
We can choose a small $\theta$ which is fixed from now on, such that
$|I_0|\leq \la_0 \ln(h^{-{\fr 1 8}}),$ where $0<r<\theta$, and
$0<\la_0<1$ which only depends on $\a$ and $\beta$.

The following two Lemmas are similar to those in [8], we give them
for completeness.
\begin{lemma} Suppose that $u\geq 0$ is a solution of
equation (1.5) in ${\cal B}^-_r$ centered at $(0,0)$ and $
mes\{(x,y,t)\in {\cal B}^-_r, \, u \geq 1\} \geq \fr 1 2 mes ({\cal
B}^-_r). $ Then there exist constants $\theta$ and $h_0$, $0<\theta,
h_0<1$ which only depend on
 $\la_0$ and $\mu$, such that
$$
u(x,y,t) \geq h_0\quad \hbox{in}\quad {\cal B}^-_{\theta r}. \leqno(3.17)
$$
\end{lemma}
{\it Proof:} We consider $$w=\ln^+(\fr{h}{u+h^{\fr98}}),$$ for
$0<h<1$, to be decided. By applying Lemma 3.3 to $w$, we have $$
-\!\!\!\!\!\!\int_{{\cal B}^-_{\theta r}}( w-I_0)_+^2 \leq
C\fr{\theta r^2}{|{\cal B}^-_{\theta r}|} \int_{{\cal
B}^-_r }|\ptl_xw|^2.
$$
Let $\tilde{u}={\fr u h}$, then $\tilde{u}$ satisfies the conditions
of Lemma 3.1. We can get similar estimates as (3.2), (3.3), (3.4)
and (3.5), hence we have
$$
\ba{lllllll} && C\fr{\theta r^2}{|{\cal B}^-_{\theta r}|}
\int_{{\cal B}^-_{ r }}|\ptl_x w|^2\\ \\
& &\leq C(\mu)\fr{\theta r^2}{|{\cal B}^-_{\theta
r}|}[\fr{C(\mu)}{\beta^4(1-\beta)^2} +\fr 45\ln(h^{-\fr 18})]
mes(K_{\beta^{-1} r})\\ \\
& &\leq C(\theta,\mu,\beta) \ln(h^{-\fr 18}),
 \ea\leqno(3.18)
$$
where $\theta$ has been chosen. By $L^{\infty}$ estimate, there
exists a constant, still denoted by $\theta$, such that for $z \in
{\cal B}^-_{\theta r}$,
$$
w-I_0\leq C(\mu,\beta) (\ln(h^{-\fr 18}))^{\fr 12} .\leqno(3.19)
$$
Therefore we may choose $h_0$ small enough, so that
$$
C(\mu,\beta) (\ln (\fr {1}{h_0^{\fr 18}}))^{\fr 12}\leq  \ln (\fr
{1}{2h_0^{\fr 18}})-\lambda_0\ln (\fr {1}{h_0^{\fr 18}}).
$$
Then (3.13) and (3.19) derive
$$
\max_{{\cal B}^-_{\theta r}}\fr{h_0}{u+h_0^{\fr 98}}\leq \fr
{1}{2h_0^{\fr 18}},$$ which implies $\min_{{\cal B}^-_{\theta
r}}u\geq h_0^{\fr 98}$, then we finished the proof of this Lemma.
\begin{lemma}
Suppose that $u$ is a weak solution of equation (1.5) in
${\cal B}^-_r$, then exists constant $h_2$, $0<h_2<1$, such that
$$ Osc_{{\cal B}^-_{\theta r}}u\leq
h_2 \,Osc_{{\cal B}^-_{r}}u, \leqno(3.20)$$ where $\theta$ is given
in Lemma 3.5.
\end{lemma}
{\it Proof:} We may assume that
$M=\max_{{\cal B}^-_{r}}(+u)=\max_{{\cal B}^-_{r}}(-u)$, otherwise
we replace $u$ by $u-c$, since $u$ is bounded locally. Then either
$1+\fr u M$ or $1-\fr u M$ satisfies the assumption of Lemma 3.5,
and we suppose $1+\fr u M$ does, thus Lemma 3.5 implies that there
exists $h_0>0$ such that
$$\inf_{{\cal B}^-_{\theta r}}(1+\fr u M)\geq h_0,\, $$ that is $u\geq
M(h_0-1)$, then
$$
Osc_{{\cal B}^-_{\theta r}}u\leq M-M(h_0-1)\leq
(1-\fr{h_0}{2})Osc_{{\cal B}^-_{r}}u,
$$
where we can let $h_2=(1-\fr{h_0}{2})$.

{\bf Proof of Theorem 2.2.} By the standard regularity arguments,
for example, see Chapter 8 in [4],
 we can obtain the result near point $(0,0)$.
By the left invariant translation group action,
 we know that $u$ is $C^{\a}$ in the interior.

\mysection{Proof of Main Theorem}  By Theorem 2.2 and a variable
transformation (see [9]), we can prove Theorem 1.1.
\\
{ \textbf{Proof}}: Since $\ptl_xb(x,y,t)\neq 0$, let $b=\xi$,
$y=\eta,$ and $t=\tau$, then
$$
\fr{\ptl}{\ptl x}=\fr{\ptl b}{\ptl x}\fr{\ptl}{\ptl \xi},\quad
\fr{\ptl}{\ptl y}=\fr{\ptl}{\ptl \eta}+\fr{\ptl b}{\ptl
y}\fr{\ptl}{\ptl \xi},\quad \fr{\ptl}{\ptl t}=\fr{\ptl}{\ptl
\tau}+\fr{\ptl b}{\ptl t}\fr{\ptl}{\ptl \xi},\quad$$
and the
equations (1.1) can be written as
$$
\displaystyle \fr{\ptl b}{\ptl x}\fr{\ptl}{\ptl \xi}(a \fr{\ptl
b}{\ptl x}\fr{\ptl}{\ptl \xi} \,u )+(b_0\fr{\ptl b}{\ptl
x}+\xi\fr{\ptl b}{\ptl y}-\fr{\ptl b}{\ptl t})\fr{\ptl}{\ptl \xi}
u+\xi\ptl_{\eta} u-{\ptl_{\tau} \,u}=0.\leqno(4.1)
$$
From the implicit function theorem, we know
$(x,y,t):\rightarrow(\xi,\eta,\tau)$ is a $C^2$ diffeomorphism, and
$\fr{\ptl b}{\ptl x}=(\fr{\ptl x(\xi,\eta,\tau)}{\ptl \xi})^{-1},$
hence the above equation can attain
$$
\displaystyle \fr{\ptl}{\ptl \xi}(\fr{\ptl b}{\ptl x}a \fr{\ptl
b}{\ptl x}\fr{\ptl}{\ptl \xi} \,u )+[a \fr{\ptl^2
x}{\ptl\xi^2}(\fr{\ptl b}{\ptl x})^3+b_0\fr{\ptl b}{\ptl
x}+\xi\fr{\ptl b}{\ptl y}-\fr{\ptl b}{\ptl t}]\fr{\ptl}{\ptl \xi}
u+\xi\ptl_{\eta} u-{\ptl_{\tau} \,u}=0,\leqno(4.2)
$$
that is
$$
\displaystyle \fr{\ptl}{\ptl \xi}(\tilde{a} \fr{\ptl}{\ptl \xi} \,u
)+\tilde{b_0}\fr{\ptl}{\ptl \xi} u+\xi\ptl_{\eta} u-{\ptl_{\tau}
\,u}=0,\leqno(4.3)
$$
where $\tilde{a}=\fr{\ptl b}{\ptl x}a \fr{\ptl b}{\ptl x}$, and
$\tilde{b_0}=a \fr{\ptl^2 x}{\ptl\xi^2}(\fr{\ptl b}{\ptl
x})^3+b_0\fr{\ptl b}{\ptl x}+\xi\fr{\ptl b}{\ptl y}-\fr{\ptl b}{\ptl
t}$.

In a fixed bounded domain, $\tilde{a}\in L^{\infty}$ and
$\tilde{b_0}\in L^{\infty}$, by Theorem 2.2 the weak solution of
(4.3) is H\"older continuous.

We give an immediate corollary.  Let $x=(x_1,\dots,x_m)$, $y=(y_1,\dots,y_n)$ and $m\geq n$.
$$
\displaystyle  Lu \equiv
\sum_{i,j=1}^m\ptl_{x_i}(a_{ij}(x,y,t)\ptl_{x_j} \,u
)+\sum_{k=1}^mb_0^k(x,y,t)\ptl_{x_k}
u+\sum_{l=1}^nb_l(x,y,t)\ptl_{y_l} u-{\ptl_t \,u}=0,\leqno(4.4)
$$
and we assume:

[H.1] the coefficients $a_{ij}$, $1\leq i,j\leq m$, are real valued,
measurable functions of $(x,t)$. Moreover, $a_{ij}=a_{ji} \in
L^{\infty} ({R}^{m+n+1})$ and there exists a $\mu
>0$ such that
$$
\mu\sum_{i=1}^{m}\xi_i^2 \leq \sum_{i,j=1}^{m}
a_{ij}(x,y,t)\xi_i \xi_j \leq \fr{1}{\mu} \sum_{i=1}^{m}\xi_i^2
$$
for every $(x,y,t)\in {R}^{m+n+1}$, and $\xi \in {R}^{m}$.

[H.2] $b_0^j\in L^{\infty}(\Omega)$, $b_l\in C^2(\Omega)$, and
$|b_0^j|_{\infty},|b_l|_{C^2}\leq \fr {1}{\mu}$, where $j=1,\dots,m$
and $l=1,\dots,n$. There exists $i_1,\cdots,i_n$, such that
$\fr{\ptl(b_1,\dots,b_n)}{\ptl(x_{i_1},\dots,x_{i_n})}\neq 0$, where
${i_j}\in \{1,\dots,m\}$, $j=1,\dots,n$ and ${i_1}< \cdots< {i_n}$.

\begin{cor}
Under the assumption [H.1] and [H.2], the weak solutions of (4.4)
are H\"older continuous.
\end{cor}

\end{document}